\documentclass[12pt,singlespace,oneside]{article}
\usepackage{amsmath}
\usepackage{latexsym}

\input{tcilatex}

\begin{document}

\title{Real Version of Calculus of Complex Variable (II): Cauchy's Point of
View\thanks{%
Dedicated to Dr. Guillermo Gil, for the indirect but important influence he
had in the author's persistent search for profound explanations.} \ }
\author{Jose G. Vargas\thanks{%
PST\ Associates. 138 Promontory Rd., Columbia, SC 29209, USA.
josegvargas@earthlink.net} \ }
\date{}
\maketitle

\begin{abstract}
As was the case in a previous paper, the differential form $x+ydxdy$ plays
the role $z$ in the standard calculus of complex variable. The role of
holomorphic functions will now be played by strict harmonic differential
forms in the K\"{a}hler algebra of the real plane. These differential forms
satisfy the Cauchy-Riemann relations.

No new concept of differentiation is needed, and yet this approach parallels
standard Cauchy theory, but more simply. The power series and theorem of
residues come here at the end, unlike in the previous paper.
\end{abstract}

\section{Introduction}

In a previous paper, we continued developing corollary developments of
Stokes theorem in multiply connected regions of the real plane. Through the
use of K\"{a}hler algebra (Clifford algebra of differential forms), it was
found that one does not need the calculus of complex variable in order to
handle the integrations for which a physicist uses it. The focus was on
closed differential 1-forms, functions with power expansions being of the
essence very early in the argument. The first major result was the theorem
of residues.

We shall now give a different version of real calculus to replace Cauchy's
theory. The focus will be on Even DIifferential Forms (``edifs'', $u+vdxdy$)
in the K\"{a}hler algebra in the real plane and, more specifically, on the
strict harmonic ones. In this paper, we shall be oblivious to whether our
integrands admit power expansions. That is a late development in the present
context. It is left for interested parties to develop as in the standard
complex variable calculus but with our simpler concepts.

The role of complex variable will again be played by $z\equiv x+ydxdy$ and
by the relations

\begin{equation}
d\phi =\frac{1}{z}dy,\text{ \ \ \ \ \ \ }d\rho =\frac{\rho }{z}dx.
\end{equation}%
An even smaller amount of K\"{a}hler algebra is now needed, but some minor
concept in K\"{a}hler calculus is required \cite{K62}.

\section{Cauchy-Riemann and Cauchy-Goursat}

\subsection{Cauchy-Riemann equations}

Clifford product will be indicated by juxtaposition. To avoid equivocity, we
use the symbol $\partial $ for the operator%
\begin{equation}
\partial \equiv dx\dfrac{\partial }{\partial x}+dy\dfrac{\partial }{\partial
y},
\end{equation}%
and define%
\begin{equation}
d\alpha \equiv \partial \wedge \alpha ,\text{ \ \ }\delta \alpha \equiv
\partial \cdot \alpha .
\end{equation}%
We thus have%
\begin{equation}
\text{ }\partial \alpha =d\alpha +\delta \alpha .
\end{equation}%
$d\alpha $ is the standard exterior derivative. When $\alpha $ is a $1$%
-form, $\delta \alpha $ is the scalar divergence.

If $u$ and $v$ are differentiable, $\partial (u+vdxdy)$ exists. A
differential form that satisfies $\partial \alpha =0$ is called strict
harmonic. Denote Strict Harmonic edifs as $shedifs$. These admit
multiplicative inverses and satisfy the Cauchy-Riemann relations:%
\begin{equation}
\partial w=0\Longleftrightarrow \text{ \ \ }du=-\delta (vdxdy)\text{ \ \ }%
\Longleftrightarrow \text{ \ \ }u,_{x}=v,_{y};\;u,_{y}=-v,_{x}.
\end{equation}%
Functions of edifs, and of shedifs in particular, can be defined as is usual
in Clifford algebra, and thus as in the calculus of complex variable.
Polynomial, rational, exponential, trigonometric, hyperbolic and logarithmic
functions of shedifs are themselves shedifs.

\subsection{Valuations and Cauchy-Goursat theorem}

We replace integration $\int_{c}f(z)dz$ on a curve of the complex plane with
valuation of $f(x+ydxdy)$ on a curve $c$ of the real plane. The valuation $%
\left\langle w\right\rangle _{c}$ of an edif $w$ on a curve $c$ is defined
as the edif%
\begin{equation}
\left\langle w\right\rangle _{c}\equiv \left[ \int_{c}wdx\right] +dxdy\left[
\int_{c}wdy\right] ,
\end{equation}%
with momentary use of square brackets for emphasis. Easy calculations yield%
\begin{equation}
\left\langle w\right\rangle _{c}=\int_{c}(udx-vdy)\;+\;dxdy\int_{c}(udy+vdx).
\end{equation}%
The integrability conditions for these integrals to not depend on $c$ but
only on the end points of the curve are the Cauchy-Riemann relations. Thus
potentials 
\begin{equation}
U=\int udx-vdy,\text{ \ \ \ \ \ \ \ \ \ \ }V=\int (udy+vdx)
\end{equation}%
exist, which imply the existence of ``valuation potentials'' of shedifs,%
\begin{equation}
\left\langle w\right\rangle =U+Vdxdy.
\end{equation}%
The valuation potential is determined only up to an additive constant
shedif. It follows that the valuation of a shedif on a closed curve on those
domains is zero. This is the Cauchy-Goursat theorem for shedifs.

$W$ is a shedif since%
\begin{equation}
dU=udx-vdy,\;\;\;dV=udy+vdx,
\end{equation}%
which implies%
\begin{equation}
U,_{x}=u=V,_{y}\;\ \ \ \ \ \ \ \ \ \ \ U,_{y}=-v=V,_{x}.
\end{equation}%
The valuation plays the role played by integration in the calculus of
complex variable.

In domains that are not simply connected, we surround the poles enclosed by
closed curves $C$ with equally oriented circles $c_{i}$, all of them with
the same orientation as $C$ and containing one and only one pole each. We
then have%
\begin{equation}
\left\langle w\right\rangle _{C}=\sum_{i}\left\langle w\right\rangle
_{c_{i}}.
\end{equation}

\subsection{Rationale for the introduction of the concept of valuation}

Consider the integral with integer $n$%
\begin{equation}
\oint \frac{dz}{(z-z_{0})^{n+1}}
\end{equation}%
in the standard calculus of complex variable. It is $2\pi i$ for $n=0$, and $%
0$ otherwise. Let us rewrite it as%
\begin{equation}
\oint \frac{\frac{1}{(z-z_{0})^{n}}}{z-z_{0}}dz
\end{equation}%
for potential integral as in \cite{V48} of integrals of the form%
\begin{equation}
\oint \frac{f(z)}{z-z_{0}}dx
\end{equation}%
(Recall $dz=d(x+ydxdy)=dx$). Here also, the integral (15) has meaning. With $%
f(z)$ equal to $(z-z_{0})^{-n}$, its value is zero, also for $n=0.$ It is
not equivalent to the integral (13) in the standard calculus of complex
variable. As for the restricted Cauchy's integral formula, it is not
applicable in this case because $(z-z_{0})^{-n}$ is not a differential $2-$%
form, not even for $n=0.$

The role of (13) is now played by%
\begin{equation}
\left\langle \frac{1}{(z-z_{0})^{n+1}}\right\rangle _{C}.
\end{equation}%
As per definition (6), we have%
\begin{equation}
\left\langle w\right\rangle _{c}\equiv \int_{c}\frac{1}{(z-z_{0})^{n+1}}%
dx+dxdy\int_{c}\frac{1}{(z-z_{0})^{n+1}}dy.
\end{equation}%
The first integral here is zero. and so is the second one except for $n=0$,
in which case (17), and thus (16), is $2\pi dxdy.$ This conclusion about the
value of (17) can be obtained as a direct consequence of the theorem of
residues \cite{V48}. However, if we want to develop this approach
independently of the previous one, we look at (17) from the perspective of
the multivariable calculus. For integration on circles centered at ($%
x_{0},y_{0}$), we take into account that%
\begin{equation}
\frac{1}{(z-z_{0})^{n}}=\rho ^{-n}(\cos n\phi -dxdy\sin n\phi ),
\end{equation}%
and the stated result $2\pi dxdy$ again follows.

\section{Cauchy's formulas}

\subsection{Cauchy's integral formula}

We no longer face the restriction of the previous paper for this formula.
But, unlike what was the case in the previous paper, we cannot approach the
theorem in the very expeditiously way available through the theorem of
residues, not yet derived here.

Let $f(z)$ be a shedif on a simply connected region of the real plane. The
limit at $z=z_{0}$ of $f(z)-f(z_{0})$ then is zero since the scalar and
2-form parts, $u$ and $v$, of $f(z)$ have derivatives and are, therefore,
continuous. The edif%
\begin{equation}
\frac{f(z)-f(z_{0})}{z-z_{0}}
\end{equation}%
has a first order pole at $z_{0}.$ By virtue of that continuity, we have%
\begin{equation}
\left\langle \frac{f(z)}{z-z_{0}}\right\rangle _{C}=\left\langle \frac{%
f(z_{0})}{z-z_{0}}\right\rangle _{C},
\end{equation}%
as follows from (6). The right hand side is%
\begin{equation}
\oint \frac{f(z_{0})}{z-z_{0}}dx\;+\;dxdy\oint \frac{f(z_{0})}{z-z_{0}}dy.
\end{equation}%
It can be computed explicitly on circles centered at $z_{0}$. Since $%
f(z_{0})dx$ equals $u_{0}dx-v_{0}dy$, and $f(z_{0})dy$ equals $%
u_{0}dy+v_{0}dx$, we get%
\begin{equation}
\left\langle \frac{f(z_{0})}{z-z_{0}}\right\rangle _{C}=\oint \frac{-v_{0}}{%
z-z_{0}}dy\;+\;dxdy\oint \frac{u_{0}}{z-z_{0}}dy,
\end{equation}%
where we have used that the integrations over $dx$ vanish because they are
integrations over $d\rho $ by virtue of (1). Hence, finally, 
\begin{equation}
\left\langle \frac{f(z)}{z-z_{0}}\right\rangle _{C}=2\pi
dxdy(u_{0}+v_{0}dxdy)=2\pi dxdyf(z_{0})
\end{equation}%
This equation is also given the alternative form%
\begin{equation}
f(z_{0})=\frac{1}{2\pi dxdy}\left\langle \frac{f(z)}{z-z_{0}}\right\rangle
_{C}.
\end{equation}

As an example, let $C$ denote the circle of radius $1$ centered at the pole $%
z=0$. Since $z=\frac{\pi }{2}$ lies outside that circle, we have%
\begin{equation}
\left\langle \frac{1}{z(z-\frac{\pi }{2})}\right\rangle _{C}=\left\langle 
\frac{\frac{1}{z-\frac{\pi }{2}}}{z}\right\rangle _{C}=2\pi dxdy\frac{1}{-%
\frac{\pi }{2}}=-4\pi dxdy,
\end{equation}%
which is an integral for which there was no correspondence in the
Weierstrass approach.

\subsection{Co-valuations}

The co-valuation of a shedif, $w=u+vdxdy$ ($\partial w=0)$ is defined as $%
\partial w/\partial x$,%
\begin{equation}
\frac{\partial w}{\partial x}=u,_{x}+v,_{x}dxdy=v,_{y}dydy+v,_{x}dxdy=dvdy,
\end{equation}%
and also%
\begin{equation}
\frac{\partial w}{\partial x}\equiv u,_{x}dxdx+u,_{y}dydx=dudx.
\end{equation}%
Clearly $\partial z/\partial x=1$. From (7), eliminating the subscript
because the result is curve independent, we get 
\begin{equation}
\frac{\partial }{\partial x}\left\langle w\right\rangle
=\int_{c}(u,_{x}dx-v,_{x}dy)\;+\;dxdy\int_{c}(u,_{x}dy+v,_{x}dx).
\end{equation}%
Using the Cauchy-Riemann conditions, we further obtain%
\begin{equation}
\frac{\partial }{\partial x}\left\langle w\right\rangle =w.
\end{equation}%
Up to an arbitrary constant differential form, the equalities%
\begin{equation}
\left\langle \frac{\partial w}{\partial x}\right\rangle =w=\frac{\partial }{%
\partial x}\left\langle w\right\rangle ,
\end{equation}%
are further completed.

\subsection{Cauchy's integral formula for derivatives}

Rewrite (24) as%
\begin{equation}
f(z)=\frac{1}{2\pi dxdy}\left\langle \frac{f(\zeta )}{\zeta -z}\right\rangle
_{C}.
\end{equation}%
Clearly%
\begin{equation}
\frac{\partial f(z)}{\partial x}=\frac{1}{2\pi dxdy}\left\langle \frac{%
\partial }{\partial x}\frac{f(\zeta )}{\zeta -z}\right\rangle _{C}=\frac{1}{%
2\pi dxdy}\left\langle \frac{f(\zeta )}{(\zeta -z)^{2}}\right\rangle _{C}.
\end{equation}%
Successive application yields%
\begin{equation}
\frac{\partial ^{n}f(z)}{\partial x^{n}}=\frac{n!}{2\pi dxdy}\left\langle 
\frac{f(\zeta )}{(\zeta -z)^{n+1}}\right\rangle .
\end{equation}%
This formula is equivalent to formula (42) of the previous paper, but
without the restriction that $f(z)$ be a differential $2-$form.

As an example, we have 
\begin{equation}
\left\langle \frac{1}{\left( z^{2}+1\right) ^{2}}\right\rangle
_{C}=\left\langle \frac{\frac{1}{\left( z+dxdy\right) ^{2}}}{\left(
dxdy-z\right) ^{2}}\right\rangle _{C}.
\end{equation}%
for valuation around the pole $(0,1)$, i.e. $z_{0}=dxdy$. We identify the
values $n=1$ and $f(z)=\left( z+dxdy\right) ^{-2}$ for application of (33)%
\begin{equation}
\left\langle \frac{1}{\left( z^{2}+1\right) ^{2}}\right\rangle _{C}=2\pi dxdy%
\left[ \frac{\partial }{\partial x}\frac{1}{\left( z+dxdy\right) ^{2}}\right]
_{z=dxdy}=\frac{\pi }{2}.
\end{equation}

\section{Concluding remarks}

Once Cauchy's integral formulas for this formalism have been developed, the
arguments leading to the Laurent series and the theorem of residues, as well
as to the obtaining of residues for poles of higher order totally parallels
the standard treatment of the same subjects in the standard calculus of
complex variables. Being it trivial for the cognoscenti, we do not need to
that here.

The Weierstrass and Cauchy's points of view have their respective merits.
Once the difference in approach has been made clear, the following meshing
of the two developments seems to us to be the most appropriate. Start with
Weierstrass to obtain the theorem of residues. Then, instead of tackling
Cauchy's integral formulas, go to the Laurent series as in that paper \cite%
{V48}. Only then address the restricted Cauchy's integral formulas, the
restriction being of the essence of that approach, which is real in the most
strict sense of the word. But the restriction also constitutes motivation
for the present, more general approach, with its new concepts of shedifs,
valuation and co-valuation. At this points, researchers can decide whether
it is in their best interest or their students to consider again the Laurent
series, the theorem of residues and the obtaining of residues for poles of
higher order

On the opposite, conservative end, there will be also those who think that
it is not yet the time to teach their students this part of the real
calculus. We dear say, however, that they would be making their students a
big favor by teaching their students the theorem of residues as per the
previous paper, even if that is all that they do before proceeding to teach
the standard calculus of complex variable. It may spark their imagination
and eventually become fervent seekers of the best way of do mathematics. The
may find like this author did ---alas too late--- that the calculus of
differential forms constitutes the language of choice for the mathematical
needs of theoretical physicists.

\end{document}